\documentclass[10pt,a4paper]{amsart} 
\usepackage{graphicx,latexsym,amssymb,diagrams}
\theoremstyle{plain} 
\newtheorem{thm}{Theorem}
\newtheorem{lem}[thm]{Lemma} 
\newtheorem{prop}[thm]{Proposition} 
\newtheorem{cor}[thm]{Corollary} 
\theoremstyle{defn} 
\newtheorem{defn}[thm]{Definition} 
\theoremstyle{remark} 
\newtheorem{rem}[thm]{Remark}

\newtheorem{nota}[thm]{\textbf{Notations}}
\numberwithin{equation}{section}
\numberwithin{figure}{section}
\newcommand{\bd}{\begin{description}}  
\newcommand{\ed}{\end{description}}
\newcommand{\ba}{\begin{array}}     
\newcommand{\ea}{\end{array}}
\newcommand{\bc}{\begin{center}}    
\newcommand{\ec}{\end{center}}
\newcommand{\be}{\begin{enumerate}} 
\newcommand{\ee}{\end{enumerate}}
\newcommand{\beq}{\begin{eqnarray}} 
\newcommand{\eeq}{\end{eqnarray}}
\newcommand{\beQ}{\begin{eqnarray*}}
\newcommand{\eeQ}{\end{eqnarray*}}
\newcommand{\bi}{\begin{itemize}}   
\newcommand{\ei}{\end{itemize}}

\newcommand{\ov}{\overline}
\newcommand{\ve}{\varepsilon}
\newcommand{\we}{\wedge}
\newcommand{\Y}{\mathsf{Y}}

\newcommand{\Yu}{\mathsf{Y}[z_1;z_2;z_3]}
\newcommand{\Yd}{\mathsf{Y}[z_0;z_2;z_3]}
\newcommand{\Yt}{\mathsf{Y}[z_0+z_1;z_2;z_3]}
\newcommand{\Yq}{\mathsf{Y}[z_1;z_1;z_2]}
\newcommand{\Yc}{\mathsf{Y}[s;z_1;z_2]}
\newcommand{\Ys}{\mathsf{Y}[z_2;z_1;z_3]}
\newcommand{\Yg}{\mathsf{Y}[(0,1);(0,1);(0,1)]}
\newcommand{\Yij}{\mathsf{Y}[(e_i,0);(e_j,0);(0,1)]}
\newcommand{\Yijk}{\mathsf{Y}[(e_i,0);(e_j,0);(e_k,0)]}
\newcommand{\Yi}{\mathsf{Y}[(e_i,0);(0,1);(0,1)]}
\newarrow{Noto} -----
\newarrow{To} ---->
\newarrow{Mapsto} |--->
\newarrow{Into} >--->
\newarrow{Onto} ----{>>}
\newarrow{Iso} >---{>>}
\newarrow{Inclus} C--->
\newarrow{Equal} =====
\newarrow{Dash} {dash}{dash}{dash}{dash}{dash}
\newarrow{Dashto} {dash}{dash}{dash}{dash}{>}
\begin{document}
\title[Goussarov-Habiro theory for string links]{Goussarov-Habiro theory for string links and the Milnor-Johnson correspondence}
\author[J.B. Meilhan]{Jean-Baptiste Meilhan}
\address{Research Institute for Mathematical Sciences\\
         Kyoto University, Kitashirakawa, Sakyo-ku \\
         Kyoto 606-8502, Japan}
\email{meilhan@kurims.kyoto-u.ac.jp}
\begin{abstract} 
We study the Goussarov-Habiro finite type invariants theory for framed string links in homology balls. 
Their degree $1$ invariants are computed: they are given by Milnor's triple linking numbers, the mod $2$ reduction of the Sato-Levine invariant, 
Arf and Rochlin's $\mu$ invariant. These invariants are seen to be naturally related to invariants of homology cylinders 
through the so-called Milnor-Johnson correspondence: in particular, an analogue of the Birman-Craggs homomorphism for string links is computed. 
The relation with Vassiliev theory is studied. 
\end{abstract} 
\maketitle
\section{Motivations}\label{blabla}
In the late 90's, M. Goussarov and K. Habiro independently developed a finite type invariant theory for compact 
oriented $3$-manifolds. The theory makes use of an efficient surgical calculus 
machinery called \emph{calculus of claspers} \cite{G,GGP,habi}. In particular the \emph{$Y_k$-equivalence}, an equivalence relation 
for $3$-manifolds arising from calculus of claspers, plays an important role in the understanding of the invariants.\\ 
Though it is also well-defined for manifolds with links, this aspect of the theory remains so far almost non-existing in the literature. 
In the present paper, we study the case of framed $n$-string links in homology balls. For $n=1$, this is equivalent to studying homology 
spheres with framed knots. String links play an important role in the study of knots and links \cite{HL} and happen to have nice 
properties in the theory of claspers. 
Here, we compute explicitly the degree $1$ invariants (in the Goussarov-Habiro sense) for framed string links in homology balls, using 
some versions of classical invariants, such as Milnor numbers, Sato-Levine, Arf and Rochlin invariants. 
This is the outcome of the characterization of the $Y_2$-equivalence relation for these objects. 

String links are very closely related to \emph{homology cylinders} \cite{GL,L}. 
Homology cylinders over a compact connected oriented surface $\Sigma$ can be seen as a generalization of the Torelli group of $\Sigma$. 
G. Massuyeau and the author explicitely computed their degree $1$ invariants \cite{MM} ; they are given by the natural extensions of 
the \emph{first Johnson homomorphism} and the \emph{Birman-Craggs homomorphism}, initially defined for the Torelli group \cite{BC,J1,J2}. 
On the other hand, N. Habegger showed in \cite{habe} how homology cylinders are geometrically related to string links in homology balls, 
such that the extension of the first Johnson homomorphism agrees with Milnor's triple linking numbers. 
So the problem which naturally arises is to compute, likewise, the analogue of the Birman-Craggs homomorphism for this so-called 
Milnor-Johnson correspondence. 
Our computation of degree $1$ invariants of string links in homology balls allows us to answer this question.

Like Goussarov-Habiro theory, the Vassiliev theory for (classical) string links can be defined using claspers. 
This viewpoint allows us to compare both theories. More precisely, we can relate the computation of degree $1$ invariants of string 
links in homology balls to an analogous results obtained by the author on Vassiliev invariants \cite{moi}. 
We also consider the link case, where a similar statement exists \cite{TY}. 

The paper is organized as follows. We will begin with some necessary preliminary material on clasper theory. 
We compute in \S 3 the Goussarov-Habiro degree $1$ invariants for framed string links in homology balls. 
\S 3.3 is devoted to the proof of this result, and \S 3.2 contains a precise definition of the invariants it involves. 
In \S 4, we introduce homology cylinders and study the Milnor-Johnson correspondence. 
The last section deals with Vassiliev invariants.  
\noindent \subsubsection*{Acknowledgments}
The author is supported by a JSPS Postdoctoral Fellowship. \\
Most of this paper is based on my PhD thesis. 
It is a pleasure to thank my advisor Nathan Habegger for many helpful conversations and comments. 
I also thank Kazuo Habiro and Gw\'ena\"el Massuyeau for useful discussions.  
\section{Preliminaries} \label{pre}
Throughout this paper, all $3$-manifolds will be supposed to be compact, connected and oriented. 
\subsection{A brief review of the Goussarov-Habiro theory} \label{GH}
Let us briefly recall from \cite{habi,GGP,G} the basic notions of clasper theory for $3$-manifolds with links. 
\begin{defn}
Let $\gamma$ be a $n$-component link in a $3$-manifold $M$. 
A \emph{clasper} $G$ for $(M,\gamma)$ is the embedding 
$$ G: F\rTo M $$
of a surface $F$ which is the thickening of a (non-necessarily connected) uni-trivalent graph having a copy of $S^1$ attached 
to each of its univalent vertices.  $G$ is disjoint from the link $\gamma$.\\
The (thickened) circles are called the \emph{leaves} of $G$, the trivalent vertices are called the \emph{nodes} of $G$ and we still 
call the thickened edges of the graph the \emph{edges} of $G$ .
\end{defn}
\noindent In particular, a \emph{tree clasper} is a connected clasper obtained from the thickening of a simply connected unitrivalent graph (with circles attached).\\
The \emph{degree} of a clasper $G$ is the minimal number of nodes of its connected components. 

\begin{figure}[!h]
\begin{center}
\includegraphics{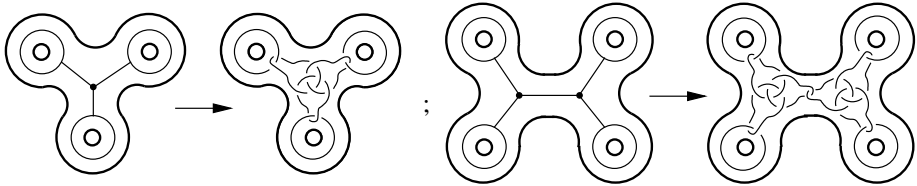}
\caption{A degree $1$ and a degree $2$ clasper and the associated framed links in their regular neighbourhoods.} \label{yh}
\end{center}
\end{figure}

A clasper $G$ for $(M,\gamma)$ is the instruction for a modification on this pair.  
There is indeed a precise procedure to construct, in a regular neighbourhood $N(G)$ 
of the clasper, an associated framed link $L_G$.  The \emph{surgery along the clasper $G$} is defined to be surgery along $L_G$. 
Though the procedure for the construction of $L_G$ will not be explained here, it is well illustrated by the two examples of Figure \ref{yh}.\footnote{
Here and throughout this paper, blackboard framing convention is used.} \\

\noindent We respectively call these two particular types of claspers \emph{$Y$-graphs} and \emph{$H$-graphs}. \\
We denote by $(M,\gamma)_G=(M_G,\gamma_G)$ the result of a surgery move on $(M,\gamma)$ along a clasper $G$:
\begin{displaymath}
\left\{ \begin{array}{ll}
\cdot & M_G = \left( M\setminus int\left( N(G_0)\right) \right) \cup_{\partial} N(G_0)_{L_G}, \\
\cdot & \textrm{$\gamma_G$ is the link in $M_G$ defined by $\gamma\subset M\setminus int\left(N(G_0)\right)\subset M_G$.}
\end{array} \right.
\end{displaymath}
\begin{defn}
Let $k\ge 1$ be an integer, and $\gamma$ be a link in a $3$-manifold $M$.
A surgery move on $(M,\gamma)$ along a connected clasper $G$ of degree $k$ is called a \emph{$Y_k$-move}. \\
The \emph{$Y_k$-equivalence}, denoted by $\sim_{Y_k}$, is the equivalence relation on $3$-manifolds with links
generated by the $Y_k$-moves and orientation-preserving diffeomorphisms (with respect to the boundary).
\end{defn} 
\noindent Note that $Y_1$-moves originally appear in \cite{matv} under the name of \emph{Borromean surgery} (as Fig. \ref{yh} suggests). 
The next proposition outlines a couple of key facts about this equivalence relation. 
\begin{prop} \label{propYk}
\be
\item[(1)] Tree claspers do suffice to define the $Y_k$-equivalence.
\item[(2)] If $1\le k\le n$, the $Y_n$-equivalence relation implies the $Y_k$-equivalence.
\ee
\end{prop}
We conclude this section with the definition of the Goussarov-Habiro theory, based on the notion of clasper. 
Consider a link $\gamma_0$ in a $3$-manifold $M_0$, and the $Y_1$-equivalence class $\mathcal{M}_0$ of $(M_0,\gamma_0)$. 
\begin{defn} \label{deffti}
Let $A$ be an Abelian group, and $k\ge 0$ be an integer. 
A \emph{finite type invariant of degree $k$} (in the Goussarov-Habiro sense) on $\mathcal{M}_0$ is a map $f:\mathcal{M}_0\rightarrow A$ such that, 
for all $(M,\gamma)\in \mathcal{M}_0$ and all familly $F=\{G_1,...G_{k+1}\}$ of $(k+1)$ disjoint $Y$-graphs for $(M,\gamma)$, the following equality holds:
$$ \sum_{F'\subseteq F} (-1)^{|F'|} f\left((M,\gamma)_{F'}\right) = 0. $$
\end{defn}
\subsection{Vassiliev theory using claspers} \label{VASS}
Another aspect of the theory of claspers is that it allows to redefine and study Vassiliev invariants of knots and links in a 
\emph{fixed} manifold \cite{habi,G2}. Here, for simplicity, we recall the definitions for the case of knots in $S^3$. 
For more about Vassiliev invariants, see \cite{BN}.

\begin{figure}[!h]
\begin{center}
\includegraphics{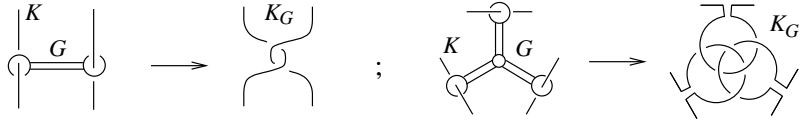}
\caption{A $C_1$-move and a $C_2$-move.} \label{exa}
\end{center}
\end{figure}

\begin{defn}
Let $K$ be a knot in $S^3$. A \emph{clasper} $G$ for $K$ is the embedding 
$$ G: F\rTo S^3 $$
of a surface $F$ which is the planar thickening of a uni-trivalent \emph{tree} (a graph without loops). 
The (thickened) 1-vertices are called the \emph{disk-leaves} of $G$, and the thickened trivalent vertices and edges of the graph are still called 
\emph{nodes} and \emph{edges} respectively.
$K$ is disjoint from $G$, except from the disk-leaves which it may intersect transversely once.\\
The \emph{$C$-degree} of a connected clasper $G$ is the number of nodes \emph{plus 1}. 
\end{defn}
Again, a clasper $G$ for $K$ is the instruction for a surgical modification: it maps $K$ to a new knot $K_G$ in $S^3$. 
Examples are given for low $C$-degrees in Fig. \ref{exa}.
\begin{defn} \label{defck}
Let $k\ge 1$ be an integer, and $K$ be a knot in $S^3$.
A surgery move on $K$ along a connected $C$-degree $k$ clasper $G$ is called a \emph{$C_k$-move}. \\
The \emph{$C_k$-equivalence}, denoted by $\sim_{C_k}$, is the equivalence relation on knots 
generated by the $C_k$-moves and isotopies.
\end{defn} 

\noindent As in Prop \ref{propYk}(2), the $C_n$-equivalence relation implies the $C_k$-equivalence if $1\le k\le n$.
\begin{rem} \label{Ckmoves}
(1). Note that a $C_1$-move is just a crossing change. 
As \cite[Fig. 2.2]{MN} shows, a $C_2$-move is equivalent to a $\Delta$-move. 
Moreover, a $C_3$-move is equivalent to a clasp-pass move (see \S \ref{ghv} for a definition) \cite{habi}.\\
(2). The $C_{k+1}$-equivalence implies the $Y_{k}$-equivalence, for all $K\ge 1$. 
More precisely, a $C_{k+1}$-move can be regarded as a special case of $Y_{k}$-move, where the leaves of the degree $k$ clasper 
are ($0$-framed) copies of the meridian of the knot. 
\end{rem}
A $C_1$-move being equivalent to a crossing change, we can reformulate the notion of Vassiliev invariant in terms of claspers.
\begin{defn} \label{defvass}
Let $A$ be an Abelian group, and $k\ge 0$ be an integer. 
An $A$-valued knot invariant $v$ is a \emph{Vassiliev invariant of degree $k$} if, 
for all knot $K$ and all familly $F=\{C_1,...C_{k+1}\}$ of $(k+1)$ disjoint $C$-degree 1 claspers for $K$, the following equality holds:
$$ \sum_{F'\subseteq F} (-1)^{|F'|} v\left(K_{F'}\right) = 0. $$
\end{defn}
\section{Goussarov-Habiro theory for string links in homology balls.} \label{GHSL}
Here and throughout the paper, unless said otherwise, by homology we mean \emph{integral} homology.  
Thus by homology ball we mean a compact oriented 
$3$-manifold whose integral homology groups are isomorphic to those of the $3$-ball.  
\subsection{String links in homology balls} \label{sl}
%
%
\subsubsection{Definition and properties}
Let $D^2$ be the standard two-dimensional disk, and $x_1,...,x_n$ be $n$ marked points in the interior of $D^2$. 
\begin{defn} \label{defsl}
An \emph{$n$-component string link} in a homology ball $M$, also called \emph{$n$-string link}, is a proper, smooth embedding 
$$ \sigma : \bigsqcup_{i=1}^n I_i \rTo M $$
of $n$ disjoint copies $I_i$ of the unit interval such that, for each $i$, the image $\sigma_i$ of $I_i$ runs from
$(x_i,0)$ to $(x_i,1)$ \emph{via} the identification $\partial M=\partial(D^2\times I)$. \\
$\sigma_i$ is called the \emph{$i^{th}$ string of $\sigma$}. It is equipped with an (upward) orientation induced by the natural orientation of $I$.\\
A \emph{framed $n$-string link} in $M$ is a string link equipped with an isotopy class of non-singular sections of its normal bundle, whose
restriction to the boundary is fixed. 
%
%
\end{defn}
\noindent We denote by $\mathcal{SL}^{hb}(n)$ the set of framed $n$-string links in homology balls, considered up to diffeomorphisms 
relative to the boundary (that is, up to diffeomorphisms whose restriction to the boundary is the identity). \\
Given two elements $(M,\sigma)$ and $(M',\sigma')$ of $\mathcal{SL}^{hb}(n)$, we can define their \emph{product} as follows. 
Denote by $M\cdot M'$ the homology ball obtained by identifying 
$\Sigma\times \{1\}\subset \partial M$ and $\Sigma\times \{0\}\subset \partial M'$.
$(M,\sigma)\cdot (M',\sigma')$ is defined by stacking $\sigma'$ over $\sigma$ in $M\cdot M'$. \\
\begin{figure}[!h]
\begin{center}
\includegraphics{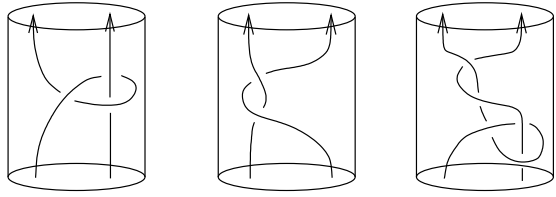}
\caption{Two $2$-string links in $D^2\times I$, and their product.} \label{slprod}
\end{center}
\end{figure} 

\noindent This product induces a monoid structure on $\mathcal{SL}^{hb}(n)$, with $(D^2\times I,1_n)$ as unit element. 
Here $1_n$ is the trivial $n$-string link.
\begin{nota} \label{notations}
Throughout this paper, the notation $1_{D^2}$ will be often used for the product $D^2\times I$. \\
$D^2_n$ will denote the $n$-punctured disk $D^2\setminus \{x_1,...,x_n\}$. 
$H := H_1(D^2_n,\mathbf{Z})$ will denote its first integral homology group, and $H_{(2)} := H_1(D^2_n,\mathbf{Z_2})$. \\
$\mathcal{B}=\{e_1,...,e_n\}$ denotes the basis of $H$ induced by the $n$ curves $h_1$, $h_2$, ... ,$h_n$ of $D^2_n$ shown in 
Fig. \ref{vatfair}.  \\
\begin{figure}[!h]
\begin{center}
\includegraphics{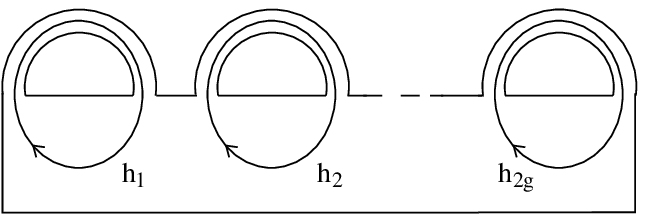}
\caption{} \label{vatfair}
\end{center}
\end{figure}
\noindent Similarly, $\mathcal{B}_{(2)}=\{\ov{e_1},...,\ov{e_n}\}$ is the associated basis of $H_{(2)}$. 
\end{nota}
Let $(M,\sigma)\in \mathcal{SL}^{hb}(n)$. 
We denote by $\hat{M}$ the homology sphere obtained by pasting a copy of $(D^2\times I)$ along its boudary, \emph{via} the 
identification $\partial M=\partial(D^2\times I)$. 
At the string links level, suitably pasting a copy of $(1_{D^2},1_n)$ along the boundary of $M$ maps $\sigma\subset M$ to a 
framed oriented $n$-component link $\hat{\sigma}\subset \hat{M}$. 
$(\hat{M},\hat{\sigma})$ is called the \emph{closure} of $(M,\sigma)$. 
In particular, for $M=1_{D^2}$, it is the usual notion of closure for $\sigma$ as defined in \cite{HL}.

\noindent Given an element $(M,\sigma)$ of $\mathcal{SL}^{hb}(n)$, let us denote by $T(\sigma)$ a tubular neighbourhood of $\sigma$. 
We denote by $M^{\sigma} := M\setminus T(\sigma)$ the exterior of the string link: the boundary of $M^{\sigma}$ is identified with 
$\partial(D^2_n\times I)$. Let $i_\epsilon$ ($\epsilon=0,1$) be the embeddings 
$$ i_{\epsilon} : D_n^2 \rTo D_n^2\times \{\epsilon\}\subset M^{\sigma}.$$
We need the following classical result of Stallings.  
\begin{thm}{\cite[Thm. 5.1]{stall}}
Let $f:A\rightarrow B$ be a map between connected CW-complexes that induces an isomorphism on the first homology groups  
and a surjective homomorphism on the second homology groups.  
Then for all $k\ge 2$, $f$ induces an isomorphism at the level of each nilpotent quotient of the fundamental group 
$$f_{k} : \frac{\pi_1(A)}{(\pi_1(A))_k}\rTo^{\simeq} \frac{\pi_1(B)}{(\pi_1(B))_k},$$
where, for any group $G$, $G_k$ is the $k^{th}$ term of its lower central series. 
\end{thm}
So by a standard Mayer-Vietoris calculation and the above theorem, the map $i_{\epsilon}$ ($\epsilon=0,1$) induces an isomorphism 
$$(i_{\epsilon})_{k} : \frac{\pi_1(D_n^2)}{(\pi_1(D_n^2))_k}=\frac{F}{F_k} 
\rTo^{\simeq} \frac{\pi_1(M^{\sigma})}{(\pi_1(M^{\sigma}))_k},$$
for each $k\ge 2$, where $F$ stands for the free group on $n$ generators.  
So any element $\sigma$ of $\mathcal{SL}^{hb}(n)$ induces an automorphism of $F/F_{k+1}$, called its \emph{$k^{th}$ Artin 
representation}, defined by $\mathcal{A}_k(\sigma)=(i_1)_{k+1}^{-1}\circ (i_0)_{k+1}$. 

\noindent Actually, $\mathcal{A}_k(\sigma)$ conjugates each generator $x_i$ of $F/F_{k+1}$ by $\lambda_i$, the \emph{$i^{th}$ longitude 
of $\sigma$ mod $F_{k+1}$}: the framing on $\sigma$ defines a curve in $M^{\sigma}$ parallel to $\sigma_i$, which determines an 
element of $\pi_1(M^{\sigma})$.  The image in $F/F_{k+1}$ of this element by $(i_1)_{k+1}^{-1}$ is $\lambda_i$. 

\noindent Denote by $\mathcal{SL}^{hb}(n)[k]:=Ker\mathcal{A}_{k}$ 
the submonoid of all $n$-string links inducing the identity on $F/F_{k+1}$. 
Note that $\mathcal{SL}^{hb}(n)=\mathcal{SL}^{hb}(n)[1]$ and that $(M,\sigma)\in \mathcal{SL}^{hb}(n)[2]$ if and only if 
$\sigma$ has null-homologous longitudes, that is, vanishing framings and linking numbers \cite{habmas}.  
\subsubsection{Goussarov-Habiro theory for framed string links in homology balls} \label{state}
Denote by $\mathcal{SL}_k^{hb}(n)$ the submonoid of all elements $(M,\sigma)\in \mathcal{SL}^{hb}(n)$ which are $Y_k$-equivalent to 
$(1_{D^2},1_n)$. There is a descending filtration of monoids
$$ \mathcal{SL}^{hb}(n)\supset \mathcal{SL}_1^{hb}(n)\supset \mathcal{SL}_2^{hb}(n)\supset ... $$
and for all $k\ge 1$, the quotient 
$$ \ov{\mathcal{SL}}_k^{hb}(n) := \mathcal{SL}_k^{hb}(n)/Y_{k+1} $$
is an Abelian group (this follows from standard calculus of claspers). 
This section is devoted to the study of the case $k=1$. First, we identify the monoid  $\mathcal{SL}_1^{hb}(n)$.
\begin{prop} \label{sl1}
The elements of $\mathcal{SL}_1^{hb}(n)$ are those $n$-string links in homology balls with vanishing framings and linking numbers:  
$$\mathcal{SL}_1^{hb}(n)=\mathcal{SL}^{hb}(n)[2].$$ 
\end{prop}
\noindent (the proof is postponed to the end of this section). 
The next result characterizes the degree $1$ Goussarov-Habiro finite type invariants 
for string links in homology balls.
\begin{thm} \label{corY2}
Let $(M,\sigma)$ and $(M',\sigma')$ be two $n$-string links in homology balls with vanishing framings and linking numbers (\textit{i.e.} two 
elements of $\mathcal{SL}_1^{hb}(n)$). The following assertions are equivalent:
\begin{enumerate}
\item[(a)] $(M,\sigma)$ and $(M',\sigma')$ are  $Y_2$-equivalent;
\item[(b)] $(M,\sigma)$ and $(M',\sigma')$ are not distinguished by degree $1$ Goussarov-Habiro finite type invariants;
\item[(c)] $(M,\sigma)$ and $(M',\sigma')$ are not distinguished by Milnor's triple linking numbers, nor the mod 2 reduction of 
the Sato-Levine invariant, the Arf invariant and Rochlin's $\mu$-invariant.
\end{enumerate}
\end{thm} 
See \S \ref{invsl} for the definitions of the above-mentioned invariants. 
\begin{rem} \label{remthm}
When considering higher degrees, the implication (a) $\Rightarrow$ (b) remains true (as well as for knots and links in homology spheres). 
The converse implication is also true when $n=1$, that is for knots in homology spheres (see \cite{habi}), and it is conjectural for 
string links with $n>1$ components. \\
This conjecture is to be compared with \cite[Conj. 6.13]{habi}, for Vassiliev invariants of (classical) string links (see also 
\S \ref{ghv}). 
\end{rem}
The proof of the theorem is given in \S \ref{y2slbh}. It consists in computing the Abelian group $\ov{\mathcal{SL}}_1^{hb}(n)$, 
in a \emph{graphical} way. 
More precisely, we will define in \S \ref{a1p} a $\Y$-shaped diagrams space $\mathcal{A}_1(P_n)$ and a surjective 
surgery map $\mathcal{A}_1(P_n)  \rTo^{\varphi_1}  \ov{\mathcal{SL}}_1^{hb}(n)$. 
We will see that $\psi$ turns out to be an isomorphism, with inverse induced by the invariants listed in Thm. \ref{corY2}.
\subsubsection{$Y_1$-equivalence for string links: proof of Proposition \ref{sl1}} \label{A.3}
We first prove the inclusion $\mathcal{SL}_1^{hb}(n)\subset  \mathcal{SL}^{hb}(n)[2]$: any element of $\mathcal{SL}^{hb}(n)$ obtained from $(1_{D^2},1_n)$ 
by a finite sequence of $Y_1$-moves has null homologous longitudes. 
It suffices to show that, if $(M_2,\sigma_2)$ is obtained from $(M_1,\sigma_1)\in \mathcal{SL}^{hb}(n)$ by surgery along a degree $1$ clasper $G$, 
these elements have homologous longitudes. 
Denote by $M_i^{\sigma_i}$ the exterior of the string links ($i=1,2$). We have 
$$ M_2^{\sigma_2} \cong \left(M_1^{\sigma_1}\right)\setminus int(N(G)) \cup_{j|_{\partial}\circ h} (H_3), $$
where $j : H_3\rInclus 1_{D^2}\setminus 1_n$ is the embedding of a genus $3$ handlebody onto a regular neighbourhood $N(G)$ of $G$, and where $h$ is an 
element of the Torelli group of $\Sigma_3=\partial H_3$ -- see \cite[Lem. 1]{Mass} for an explicit description of this diffeomorphism. 
$h$ induces the identity on $\pi_1(\Sigma_3)/\pi_1(\Sigma_3)_2$: it follows (by a Van Kampen type argument) that there is an isomorphism 
$$ \frac{\pi_1(M_1^{\sigma})}{\left(\pi_1(M_1^{\sigma})\right)_2} \rTo^{\simeq } \frac{\pi_1(M_2^{\sigma'})}{\left(\pi_1(M_2^{\sigma'})\right)_2},$$
which is compatible with the maps $i_{\ve}$ ; $\ve=0,1$. The assertion follows. 

The other inclusion is essentially due to N. Habegger \cite{habe}. First, recall that every homology sphere is $Y_1$-equivalent to the $3$-sphere $S^3$ 
\cite{matv,habi} ; likewise every homology ball is $Y_1$-equivalent to $B^3\cong D^2\times I$. 
So it suffices to show that a framed string link $\sigma$ in $D^2\times I$ whose framings and linking numbers are all zero is $Y_1$-equivalent
to $(1_{D^2},1_n)$. By a sequence of connected sums on $\sigma$ with copies of the $0$-framed Borromean link, we can furthermore suppose that 
all Milnor's triple linking numbers are zero: such connected sums are nothing else but $Y_1$-moves (each leaf of the clasper being a meridian of 
the string on which connected sum is performed). 
By \cite[Thm. D]{levi}, $\sigma$ is thus \emph{surgery equivalent} to the trivial string link, that is, 
$\sigma$ is obtained from $1_n$ by a sequence of surgeries on trivial $(\pm 1)$-framed knots $K_i$ in the exterior of $\sigma$, these knots having vanishing 
linking numbers with $\sigma$. The union $\cup_i K_i$ is a $(\pm 1)$-framed boundary link: surgery on such a link is known to be equivalent 
to a sequence of $Y_1$-surgeries \cite[Cor. 6.2]{habe}. 
\subsection{Classical invariants for string links in homology balls} \label{invsl}
\subsubsection{Rochlin's $\mu$-invariant} \label{roch}
Let $M$ be a closed $3$-manifold endowed with a spin structure $s$, and let $(W,S)$ be a compact spin $4$-manifold spin-bounded by $(M,s)$ (that is, 
$\partial W=M$ and $S$ coincides with $s$ on $M$). Then, the \textit{modulo 16} signature $\sigma(W)$ of $W$ 
is a well-defined closed spin $3$-manifolds invariant $R(M,s)$, called the \emph{Rochlin invariant} of $M$. 
In the case of homology spheres, there is a unique spin structure $s_0$, and $R(M,s_0)$ is divisible by $8$: 
$$ \mu(M) := \frac{R(M,s_0)}{8} \in \mathbf{Z}_{2} $$
is an invariant of homology spheres called \emph{Rochlin's $\mu$-invariant}. 

For elements $(M,\sigma)$ of $\mathcal{SL}^{hb}(n)$, we set 
$$ R(M,\sigma) := \mu(\hat{M}), $$ 
\noindent where the homology sphere $\hat{M}$ is the closure of $M$ as defined in \S \ref{sl}. 
The following result of G. Massuyeau implies that the restriction of $R$ to $\mathcal{SL}^{hb}_1(n)$ factors to a homomorphism of Abelian groups  
\begin{diagram}
R : & \ov{\mathcal{SL}}_1^{hb}(n) & \rTo & \mathbf{Z}_2.
\end{diagram}

\begin{prop}\cite[Cor. 1]{Mass} \label{propmu}
Rochlin's invariant is a degree $1$ finite type invariant (in the Goussarov-Habiro sense) of integral homology spheres. 
\end{prop}
\subsubsection{Milnor Invariants} \label{miln}
Let $\sigma$ be an $n$-string link in a homology ball $M$. Recall from \S \ref{sl} that $F$ is the free group on $n$ generators, 
and that $F_k$ is the $k^{th}$ term of its lower central series. 
Recall also that $\lambda_i\in F/F_{k+1}$ denotes the $i^{th}$ longitude of $\sigma$ mod $F_{k+1}$.

\noindent Denote by $P(n)$ the ring of power series in the non-commutative variables $X_1,...,X_n$. 
The \emph{Magnus expansion} \cite{MKS} $F\rTo P(n)$ is a group homomorphism which maps each generator $x_i$ of $F$ to $1+X_i$. 
\begin{defn}
The \emph{Milnor's $\mu$-invariant of length $l$}, $\mu_{i_1...i_l}$ of $\sigma$ is the coefficient of the monomial $X_{i_1}...X_{i_{l-1}}$ 
in the Magnus expansion of the longitude $\lambda_{i_l}\in F/F_k$ for a certain $k\ge l$.
\end{defn}
\noindent For example, Milnor's invariants of length $2$ are just the linking numbers. Here, we deal with  Milnor's invariants of length $3$, 
also called \emph{Milnor's triple linking number}. The following proposition-definition follows from Lemma \ref{lem:Y2sl} below.
%
\begin{prop} \label{propmu3}
For all $i<j<k\in \{1,...,n\}$, there is a well-defined homomorphism of Abelian groups
$$ \ov{\mathcal{SL}}_{1}^{hb}(n)\rTo^{\mu_{ijk}} \mathbf{Z} $$ 
induced by Milnor's triple linking number. 
\end{prop}
\begin{rem}
In general, Milnor's triple linking numbers are not additive on $\mathcal{SL}(n)$. The homomorphism defect is given by linking numbers, so it 
vanishes for elements of $\mathcal{SL}^{hb}_1(n)$.
\end{rem}
\begin{lem} \label{lem:Y2sl} 
Let $(M,\sigma)$ be a framed string link in a homology ball. Let also $G$ be a degree $2$ clasper in $M$ 
disjoint from $\sigma$ and let $(M_G, \sigma_G)$ be the result of the surgery along $G$.
Then, there exists an isomorphism
\begin{diagram}
\frac{\pi_1(M^{\sigma})}{(\pi_1(M^{\sigma}))_3} & \rTo^{\simeq} & 
\frac{\pi_1(M^{\sigma}_G)}{\left(\pi_1(M^{\sigma}_G)\right)_3} 
\end{diagram}
compatible with the inclusions $i_{\ve}$ ; $\ve=0,1$.
\end{lem}
\emph{Proof: } The reader is refered to the proof of \cite[Lem. 3.13]{MM}. $\square$ 
\subsubsection{The Arf Invariant} \label{ARF}
\begin{figure}[!h] 
\begin{center} 
\includegraphics{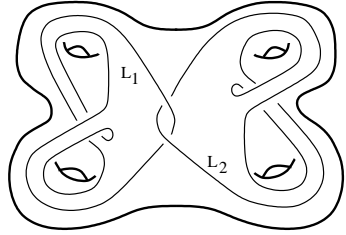}
\caption{The $2$-component link $L$.} \label{fig:link} 
\end{center} 
\end{figure}

Let $K$ be a knot in a homology sphere $M$, and $S$ be a Seifert surface for $K$ of genus $g$. 
Denote by $\centerdot$ the \emph{mod $2$ reduction} of the homological intersection form on $H_1(S,\mathbf{Z}_2)$. 
Let $\delta_2: H_1(S,\mathbf{Z}_2)\rTo \mathbf{Z}_2$ be the map defined by 
$$  \delta_2(\alpha) = lk(\alpha,\alpha^+) (mod\textrm{ 2}), $$
where $\alpha^+$ is a parallel copy of $\alpha$ in the positive normal sense of $S$ (for a fixed orientation of $M$). 
$\delta_2$ is a quadratic form with $\centerdot$ as associated bilinear form: the \emph{Arf invariant of the knot $K$} 
\cite{rob} is the Arf invariant of $\delta_2$, that is, for a given symplectic basis $\{a_1,b_1,...,a_g,b_g\}$ for $\centerdot$ 
$$ Arf(K) = Arf(\delta_2) = \sum_{i=1}^g \delta_2(a_i) \delta_2(b_i). $$
\begin{rem} \label{remdef}
The fact that the Arf invariant is still well-defined for knots in homology spheres essentially follows from 
the following fact (see for example \cite{GT} for a proof): 
two Seifert surfaces $S_0$ and $S_1$ for a knot $K$ in an homology sphere $M$ are related by a sequence of 
isotopies, additions and removals of tubes $S^1\times I$. 
Indeed, as we will see in the proof of Prop. \ref{proparf}, such tubes do not contribute to the Arf invariant.
\end{rem}
For elements of $\mathcal{SL}_1^{hb}(n)$, the Arf invariant is defined in the obvious way: for an integer $1\le i\le n$, denote by $a_i(M,\sigma)$ the 
Arf invariant of $\hat{\sigma_i}$, the $i^{th}$ component of the link $\hat{\sigma}\in \hat{M}$.  
We clearly have the following proposition-definition:
\begin{prop}
For any integer $1\le i\le n$, the map $a_i : \mathcal{SL}_1^{hb}(n)\rTo \mathbf{Z}_2$ is a homorphism of monoids, called the 
\emph{$i^{th}$ Arf invariant of $(M,\sigma)$}. 
\end{prop}
Further, this invariant happens to behave well under a $Y_2$-move. 
\begin{prop} \label{proparf}
The Arf invariant of knots in homology spheres is invariant under a $Y_2$-move. \\
As a consequence, for any $1\le i\le n$, the $i^{th}$ Arf invariant of string links in homology balls factors through a homomorphism 
of Abelian groups  
$$ a_i : \ov{\mathcal{SL}}_1^{hb}(n)\rTo  \mathbf{Z}_2. $$
\end{prop}
\emph{Proof: }
Let $K$ be a knot in a homology sphere $M$, and let $S$ be a Seifert surface for $K$. 
Let $G$ be a degree $2$ clasper for $(M,K)$ ; thanks to Prop. \ref{propYk}(2), we can suppose that $G$ is a $H$-graph. 
It suffices to show that
$ Arf(M,K)=Arf(M_G,K_G)\in \mathbf{Z}_2. $\\
Denote by $N$ a regular neighbourhood of $G$, which is a genus $4$ handlebody. 
The $10$-component surgery link associated to $G$, depicted in Figure \ref{yh}, is Kirby-equivalent to the 2-component link $L$ 
depicted in Fig. \ref{fig:link}.  This can be checked by using moves 2, 9 and 1 of \cite[Prop. 2.7]{habi} (see also \cite[pp. 254]{L}).

\begin{figure}[!h]
\begin{center}
\includegraphics{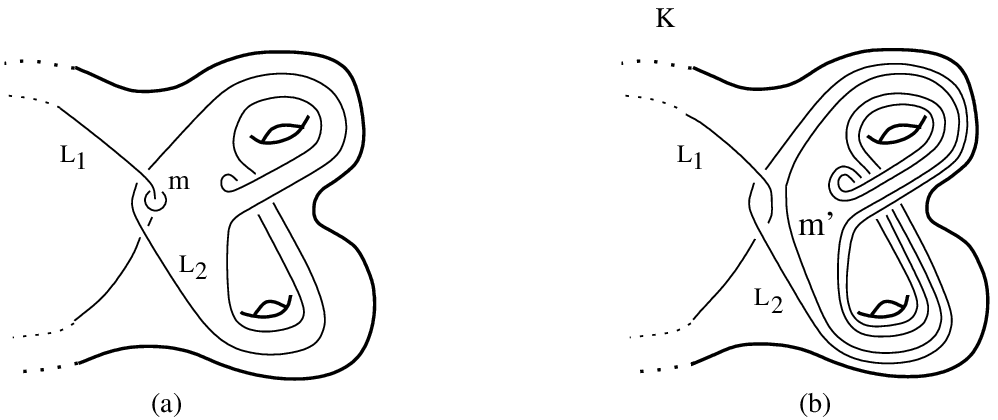}
\caption{} \label{tub}
\end{center}
\end{figure}

$K$ being disjoint from $G$, we can suppose that it is also disjoint from $N$. But $S$ may intersect $N$ and the knot $K$. 
We construct a new Seifert surface $S'$ for $L$, satisfying $S'\cap L=\emptyset$, by adding tubes $S^1\times I$ to $S$ in $N$:  
these tubes are portions of (parallel copies of) a tubular neighbourhood of the link $L$.  
The general procedure for constructing $S'$ is explained in Appendix \ref{appendix}. 
$S'$ can be seen in $M\setminus L$, and thus in the surgered manifold $M_G$.\\
Now observe that such an addition of tube doesn't affect the Arf invariant of $K$: 
if we denote by ($m,l$) a meridian/longitude pair for this tube, we have indeed $\delta_2(m)=0$, such a meridian $m$ having vanishing 
self-linking. \\
We must also show that this pair does not contribute to the Arf invariant of $(K)_G$. 
In other words, if we denote by $(m',l')$ the image of $(m,l)$ after surgery on $L$, we must show that $\delta_2(m')\delta_2(l')=0$.
Observe that the meridian $m$ can be isotoped in a small ball $B$ of $N$ where the crossing between $L_1$ and $L_2$ occurs - see Fig. 
\ref{tub}(a). Thus, surgery on $L$ sends $m$ to a curve $m'$, which is a parallel copy of $L_2$ outside of $B$, as shown in Fig. 
\ref{tub}(b): we have $\delta_2(m')=lk(m',(m')^+)=0$.\\
$\square$ \\
\subsubsection{The Sato-Levine invariant} \label{SL}

\begin{figure}[!h]
\bc
\includegraphics{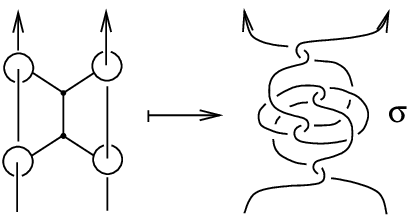}
\caption{} \label{metr}
\ec
\end{figure}

\begin{figure}[!h]
\begin{center}
\includegraphics{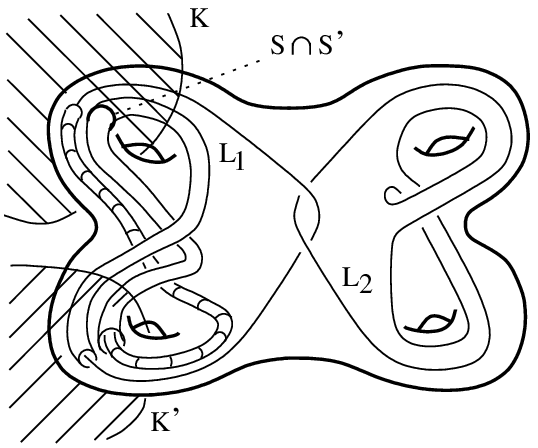}
\caption{} \label{tub2}
\end{center}
\end{figure}

Let $L=L_1\cup L_2$ be a $2$-component oriented link such that $lk(L_1,L_2)=0$. The components of $L$ bound some Seifert surfaces 
$S_1$ and $S_2$ such that $L_1\cap S_2=L_2\cap S_1=\emptyset$. $S_1$ and $S_2$ intersect along circles 
$S_1\cap S_2=C_1\cup ...\cup C_n=C$. The self-linking of $C$ relative to any of both surfaces is called the \emph{Sato-Levine invariant} 
of $L$ \cite{sato}: 
$$ \beta(L) = lk(C,C^+). $$
\noindent The fact that $\beta$ is still well-defined for links in homology spheres is again a concequence of the fact recalled in 
Rem. \ref{remdef}. Indeed, if we add a tube $t$ to (say) $S_1$, it will only intersect $S_2$ along copies of a meridian of $t$ (up to 
isotopy): such a meridian has vanishing self-linking number and links no other component of $S_1\cap S_2$.\\
The Sato-Levine invariant can also be defined for elements $(M,\sigma)$ of $\mathcal{SL}_1^{hb}(n)$. For any pair 
of integers $(i,j)$ such that $1\le i<j\le n$, we denote by $\beta_{ij}(M,\sigma)$ the Sato-Levine invariant of the 
$2$-component link of $\hat{M}$ obtained by closing the $i^{th}$ and $j^{th}$ components of $\sigma$: 
$\beta_{ij}(M,\sigma) :=\beta(\hat{\sigma_i}\cup \hat{\sigma_j})$. \\
Note that this makes sense by Prop. \ref{sl1}, as elements of $\mathcal{SL}_1^{hb}(n)$ 
have vanishing linking numbers. Moreover, $\beta_{ij}$ is additive. 
\begin{prop}
$\forall$ $1\le i<j\le n$, the map $\beta_{ij} : \mathcal{SL}_1^{hb}(n)\rTo \mathbf{Z}$ is a homomorphism of monoids,  
called the \emph{Sato-Levine invariant $\beta_{ij}$}.
\end{prop}
\noindent Note that the Sato-Levine invariant is not invariant under $Y_2$-moves: for example, it takes value $2$ on the string link $\sigma$
depicted below, obtained by surgery on $(1_{D^2},1_n)$ along a $H$ graph whose leaves are meridians of $1_n$ as depicted in Fig. \ref{metr}. 
But it turns out that it is the case for its mod $2$ reduction.
\begin{prop} \label{propsl}
The \emph{mod 2 reduction} of the Sato-Levine invariant of links in homology spheres is invariant under a $Y_2$-move. \\
In particular, for any $1\le i<j\le n$, the Sato-Levine invariant $\beta_{i,j}$ of string links in homology balls factors through 
a homomorphism of Abelian groups  
$$ \beta^{(2)}_{ij} : \ov{\mathcal{SL}}_1^{hb}(n)\rTo \mathbf{Z}_2. $$
\end{prop}
\emph{Proof: }
Let $K\cup K'$ be a $2$-component oriented link with linking number $0$ in a homology sphere $M$. Let $G$ be a degree $2$ clasper for $(M,K\cup K')$ (which, as 
in the preceding proof, can be supposed to be a $H$-graph), and $N$ be a regular neighbourhood of $G$. We must show that
$$ \beta^{(2)}(M,K\cup K')=\beta^{(2)}(M_G,K_G\cup K'_G)\in \mathbf{Z}_2. $$
\noindent We denote respectively by $S$ and $S'$ a Seifert surface for $K$ and $K'$: $S\cap S'=C_1\cup...\cup C_n=C$. 
Consider in $N$ the $2$-component surgery link $L=L_1\cup L_2$ associated to $G$ depicted in Fig. \ref{fig:link}. 
$K$ and $K'$ are supposed to be disjoint from $N$, but $S$ and $S'$ may intersect $N$ (and thus $L$).\\
When $S$ (resp. $S'$) intersects $L$, we add some tubes to built a new Seifert surface for $K$ (resp. $K'$), which is disjoint 
from $L$. The procedure for such an addition of tube is the same as the procedure explained in Appendix \ref{appendix} for a knot. 
We denote by $\tilde{C}$ the set of elements of $S\cap S'$ which are possibly created (in $N$) under this addition of tube: 
$S\cap S'=C\cup \tilde{C}$. 
A simple example of such a situation is given in Figure \ref{tub2}.

\noindent Clearly, $\tilde{C}$ is a finite number of copies of small meridians of $L_1$ and $L_2$. We clearly have 
$lk(\tilde{C},C^+)=lk(C,\tilde{C}^+)=lk(\tilde{C},\tilde{C}^+)=0$. 
It remains to prove that, \emph{after surgery} along $L$, the elements of $\tilde{C}\subset S\cap S'$ do also not contribute to 
$\beta^{(2)}(K_G\cup K'_G)$. 

$\cdot$ Suppose that $\tilde{C}=\{m\}$, where $m$ is a meridian of any of both components. 
Denote by $c$ its image after surgery on $G$: as seen in the proof of Prop. \ref{proparf}, we have $lk(c,c^+)=0$. 

\begin{figure}[!h]
\begin{center}
\includegraphics{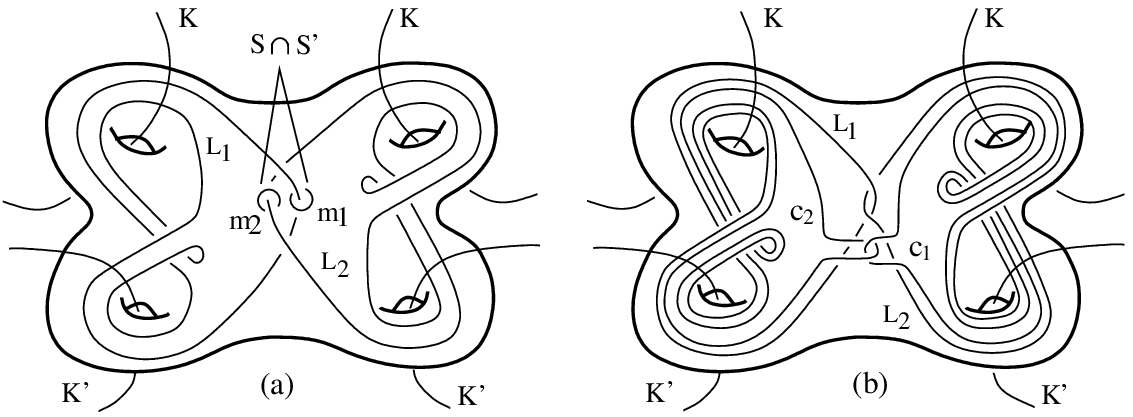}
\caption{} \label{tub3}
\end{center}
\end{figure}

$\cdot$ Now, consider the case $\tilde{C}=\{m_1,m_2\}$, a pair of meridians of $L_1$ and $L_2$. 
An example is given by the situation of Fig. \ref{tub3}(a). 

\noindent Again, 
surgery on $G$ sends $(m_1,m_2)$ to a pair of curves $(c_1,c_2)$, which are parallel copies of $L_1$ and 
$L_2$ outside of a ball of $N$ where the crossing between $L_1$ and $L_2$ occurs - see Fig. \ref{tub3}(b). Thus, $c_1$ and $c_2$ satisfy 
\beQ
lk\left(c_1\cup c_2,(c_1\cup c_2)^+\right) & = & lk(c_1,c_1^+)+lk(c_2,c_1^+)+lk(c_1,c_2^+)+lk(c_2,c_2^+)\\
 & = & 2.lk(c_1,c_2^+)=\pm 2.
\eeQ
It follows that, in these two particular cases, the mod 2 reduction of $\beta$ remains unchanged. 
The general case, where $\tilde{C}$ consists in several copies of $m_1$ and $m_2$, is proven the same way.
$\square$ \\
\begin{rem}
Note that a (less direct) proof of Prop. \ref{propsl} can be given using a formula of K. Murasugi that expresses the 
modulo $2$ reduction of the Sato-Levine invariant of a link in terms of its Arf invariants \cite{murasugi}. Indeed, 
the Arf invariant of a link can be expressed as the Arf invariant of a knot \emph{related} to $L$, that is (roughly) obtained by 
performing a connected sum of its components along some band \cite{rob}. The result then follows from Prop. \ref{proparf}. 
\end{rem}
\subsection{Degree $1$ invariants for string links: proof of Theorem \ref{corY2}} \label{y2slbh}
As announced in \S \ref{state}, the proof of Theorem \ref{corY2} consists in computing the Abelian group $\ov{\mathcal{SL}}_1^{hb}(n)$. 
This computation goes in two steps. First we will construct a combinatorial upper bound, by defining a surjective homomorphism 
$\varphi_1: \mathcal{A}_1(P_n)  \rTo \ov{\mathcal{SL}}_1^{hb}(n)$, where $\mathcal{A}_1(P_n)$ is a space of diagram. Second, 
we will show that $\psi$ is actually an isomorphism, with inverse given by the invariants listed in Thm. \ref{corY2}.\\ 
The development of the proof, and the objects it involves, are similar to those used in the proof of \cite[Thm 1.4]{MM}. We will 
recall and use several material and facts presented in the latter, to which the reader is refered for more details. 
\subsubsection{Combinatorial upper bound} \label{a1p}
Let $P_n$ denote the Abelian group $H\oplus \mathbf{Z}_2$. 
%
We denote by $\mathcal{A}_1(P_n)$ the free Abelian group generated by $\Y$-shaped unitrivalent graphs, whose trivalent vertex is equipped 
with a cyclic order on the incident edges and whose univalent vertices are labelled by $P_n$, subject to the two following relations\\[0.5cm] 
\begin{tabular}{rcrcl}  
\textbf{Multilinearity} & : &  
$\Yt$ & = & $\Yd + \Yu,$ \\[0.3cm] 
\textbf{Slide} & : & $\Yq$ & = & $\Yc,$  
\end{tabular}\\[0.5cm] 
where $z_0, z_1, z_2, z_3 \in P_n$. Here, the notation $\Y[z_1,z_2,z_3]$ stands for the graph whose univalent vertices are 
colored by $z_1$, $z_2$ and $z_3\in P_n$ in accordance with the cyclic order. This notation is invariant under cyclic permutation 
of the $z_i$'s. 
\begin{rem}
Note that, as a consequence of the Multilinearity and Slide relations, the Antisymmetry relation 
$$ \Yu = -\Ys $$
holds in $\mathcal{A}_1(P_n)$ -- for example, apply the Slide relation to $\Y[z_1+ z_2;z_1+ z_2;z_3]$. 
\end{rem}
Consider the map 
$$\rho: \mathcal{A}_1(P_n)\rTo \Lambda^{3} H\oplus \Lambda^{2} H_{(2)}\oplus H_{(2)}\oplus \mathbf{Z}_2$$ 
defined on the generators of $\mathcal{A}_1(P_n)$ by:
\bc
\begin{tabular}{rcll}
$\rho(\Yijk)$ & = & $e_i\we e_j\we e_k\in \Lambda^{3} H$, \\
$\rho(\Yij)$ & = & $\ov{e_i}\we \ov{e_j}\in \Lambda^{2} H_{(2)}$, \\
$\rho(\Yi)$ & = & $\ov{e_i}\in H_{(2)}$, \\
$\rho(\Yg)$ & = & $1\in \mathbf{Z}_2$, 
\end{tabular}
\ec
where $1\leq i < j < k \leq n$, and where $(e_i)_i$ (resp. $(\ov{e_i})_i$) are the basis elements of $H$ (resp. $H_{(2)}$) 
defined in Notations \ref{notations}. \\
$\rho$ is clearly well-defined and we actually have the following lemma. 
\begin{lem} \label{ISO}
The map $\rho$ is an isomorphism. 
\end{lem}
\noindent This is proved in the same way as \cite[Lem. 4.24]{MM} (see also \cite[Lem. 6.3]{gwen}).

%
We now construct the surgery map 
$$\varphi_1:\mathcal{A}_1(P_n) \rTo \ov{\mathcal{SL}}_1^{hb}(n). $$
For each generator $\Y=\Yu$ of $\mathcal A_1(P_n)$, where $z_i:=(h_i,\ve_i)\in P_n$, we set 
$$ \varphi_1(\Y):= (D^2\times I,1_n)_{\phi(\Y)}, $$
where $\phi(\Y)$ is a degree $1$ connected clasper (a $Y$-graph) for $(D^2\times I,1_n)$ constructed from the informations contained in the diagram $\Y$: \\
For $i\in \{ 1,2,3\}$, consider an oriented simple closed curve $c_i$ in $D^2_n\times \{1 \} \subset D^2\times I$ such that 
$[c_i]=h_i\in H$, framed along the surface. 
Then push this framed curves down in the interior of $(D^2\times I)\setminus 1_n\cong (D^2_n\times I)$, by adding a $\ve_i$-twist. 
The resulting oriented framed knot is denoted by $K_i$. 
Next, pick an embedded $2$-disk $D$ in the interior of $D^2_n\times I$ 
and disjoint from the $K_i$'s, orient it in an arbitrary way, and connect it to the $K_i$'s with some bands $e_i$. These band sums have to be compatible with 
the orientations, and to be coherent with the cyclic ordering $(1,2,3)$.
\begin{prop} \label{thchir} 
Let $\Y$ be a generator of $\mathcal A_1(P_n)$. The $Y_2$-equivalence class of $(D^2\times I,1_n)_{\phi(\Y)}$ 
does not depend on the choice of $\phi(\Y)$ (obtained by the above construction). 
Hence, we have a well-defined, surjective \emph{surgery map} 
$$ \mathcal{A}_1(P_n) \rTo^{\varphi_1} \ov{\mathcal{SL}}_1^{hb}(n).$$
\end{prop} 
The proof is strictly the same as the proof of \cite[Thm. 2.11]{MM}, and essentially uses the calculus of claspers. In particular, the 
independence on the choice of $\phi$ follows from facts similar to \cite[Cor. 4.2 and 4.3, Lem. 4.4]{GGP}.
\subsubsection{Characterization of $Y_2$-equivalence for string links} \label{prthmY}
%
Set $V:=\Lambda^2H_{(2)}\oplus H_{(2)}\oplus \mathbf{Z}_2$, and let
$$\tau : \mathcal{SL}^{hb}_1(n) \rightarrow \Lambda^3H_{(2)}\oplus V$$
be defined, for any $(M,\sigma)\in \mathcal{SL}^{hb}_1(n)$, by 
$$ \tau(M,\sigma) = \sum_{1\le i<j<k\le n} \mu^{(2)}_{ijk}(M,\sigma).\ov{e_i}\we \ov{e_j}\we \ov{e_k} 
+ \sum_{1\le i<j\le n} \beta^{(2)}_{ij}(M,\sigma).\ov{e_i}\we \ov{e_j}$$
$$ + \sum_{1\le i\le n} a_{i}(M,\sigma).\ov{e_i} + R(M). $$
Here, $\mu^{(2)}_{ijk}$ denotes the mod $2$ reduction of Milnor's triple linking number $\mu_{ijk}$.\\ 
It follows from  Propositions \ref{propmu}, \ref{propmu3}, \ref{proparf} and \ref{propsl} that this well-defined map factors 
through a homomorphism of Abelian groups 
\begin{diagram}
\ov{\mathcal{SL}}^{hb}_1(n) & \rTo^{\tau} & \Lambda^3H_{(2)}\oplus V.
\end{diagram}
Denote by $T$ the composition 
$$T: \mathcal{A}_1(P_n)\rTo^{\rho} \Lambda^3H\oplus V\rTo^{-\otimes \mathbf{Z}_2} \Lambda^3H_{(2)}\oplus V.$$
\begin{lem} \label{lem:slam}
The following diagram commutes 
\begin{diagram} 
\mathcal{A}_1(P_n) & \rOnto^{\varphi_1} & \ov{\mathcal{SL}}^{hb}_1(n) \\  
& \rdOnto<{T} & \dTo>{\tau} \\ 
& & \Lambda^3H_{(2)}\oplus V. 
\end{diagram} 
\end{lem}
\emph{Proof: }
$P$ is generated by $(0,1)$ and $(e_i,0)$, $i=1,...,n$. So, thanks to the Slide relation, there are four distinct types of generators $\Y$ for 
$\mathcal{A}_1(P_n)$, listed below ($1\le i<j<k\le n$): we prove that, in these four cases, $\tau\left(\varphi_1(\Y)\right)=T(\Y)$.

\textbf{1. $\Y=\Yg$}\\
In this case, $T(\Y)=1\in \mathbf{Z}_2$. On the other hand, 
a representative for $\varphi_1(\Y)\in \ov{\mathcal{SL}}^{hb}_1(n)$ is 
$(1_{D^2},1_n)_G$, where $G$ is contained in a ball disjoint from $1_n$ and its leaves are three copies of the $(-1)$-framed unknot. 
It follows that $(1_{D^2},1_n)_G\cong (P,1_n)$, where the closure of $P$ is the Poincar\'e sphere: $R(P,1_n) = 1$. Moreover, 
$$ \mu_{rst}\left(P,1_n \right) = \beta^{(2)}_{rs}\left(P,1_n \right) = a_{s}\left(P,1_n \right) = 0,  $$
\noindent  $\forall$ $r\ne s\ne t\in \{1,...,n\}$. It follows that $ \tau(P,1_n)=1\in \mathbf{Z}_2$.

\begin{figure}[!h]
\begin{center}
\includegraphics{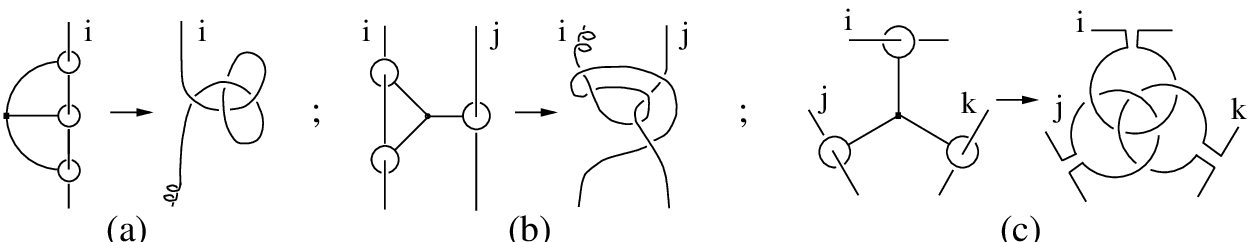}
\caption{} \label{rep}
\end{center}
\end{figure}

\textbf{2. $\Y=\Y[(e_i,0);(e_i,0);(e_i,0)]$}\\
A representative for $\varphi_1(\Y)$ is $(1_{D^2},1_n)_G$, where the three leaves 
of $G$ are small meridians of the $i^{th}$ string $(1_n)_i$ of $1_n$. Thus $(1_{D^2},1_n)_G\cong (1_{D^2},T_{i})$, where 
$T_{i}$ only differs from $1_n$ by a copy of the trefoil on the $i^{th}$ string -- see the Fig. \ref{rep}(a). We have 
$a_{r}\left(1_{D^2},T_{i} \right) = \delta_{r,i}$, and 
$$ \mu_{rst}\left(1_{D^2},T_{i} \right) = \beta^{(2)}_{rs}\left(1_{D^2},T_{i} \right) = R\left(1_{D^2},T_{i} \right) = 0
\textrm{ , $\forall$ $(r,s,t)$.} $$
It follows that $\tau\circ \varphi_1(\Y)=\ov{e_i}=T(\Y)$.

\textbf{3. $\Y=\Y[(e_i,0);(e_i,0);(e_j,0)]$}\\
A representative for $\varphi_1(\Y)$ is obtained from $(1_{D^2},1_n)$ by surgery along a $Y$-graph 
$G$ having two copies of a meridian of $(1_n)_i$ and one copy of a meridian of $(1_n)_j$ as leaves: $(1_{D^2},1_n)_G\cong 
(1_{D^2},w_{ij})$, where the $i^{th}$ and $j^{th}$ strings of $w_{ij}$ form a Whitehead link, see Fig. \ref{rep}(b). 
The Sato-Levine invariant of the Whitehead link being $1$, we obtain 
$\beta^{(2)}_{rs}\left(1_{D^2},w_{ij} \right) = \delta_{(r,s),(i,j)}$, and 
$$ \mu_{rst}\left(1_{D^2},w_{ij} \right) = a_{r}\left(1_{D^2},w_{ij} \right) = R\left(1_{D^2},w_{ij} \right) = 0
\textrm{ , $\forall$ $(r,s,t)$.} $$
It follows that $\tau\circ \varphi_1(\Y)=\ov{e_i}\we \ov{e_j}\in \Lambda^2H_{(2)}$, which coincides with $T(\Y)$.

\textbf{4. $\Y=\Y[(e_i,0);(e_j,0);(e_k,0)]$}\\
A representative for $\varphi_1(\Y)$ is $(1_{D^2},\sigma_{ijk})$, obtained from $1_n$ by performing 
a connected sum on strings $\sigma_i$, $\sigma_j$ and $\sigma_k$ with the three components of a Borromean ring, see Fig. 
\ref{rep}(c). It follows that $\mu_{abc}(\sigma_{ijk}) = 1$ for $(a,b,c)=(i,j,k)$, and $0$ otherwise. Moreover, 
$$ \beta^{(2)}_{rs}\left(1_{D^2},\sigma_{ijk} \right) = a_{r}\left(1_{D^2},\sigma_{ijk} \right) = R\left(1_{D^2},\sigma_{ijk} \right) = 0
\textrm{ , $\forall$ $(r,s)$.} $$
We thus obtain $\tau\left(\varphi_1(\Y)\right)=\ov{e_i}\we \ov{e_j}\we \ov{e_k}=T(\Y)$, 
which completes the proof.
$\square$ \\
Furthermore, we can define by Prop. \ref{propmu3} a homomorphism of Abelian groups
\begin{diagram}
\ov{\mathcal{SL}}^{hb}_1(n) & \rTo^{\mu_3} & \Lambda^3H.
\end{diagram}
by setting $\mu_3(M,\sigma) = \sum_{1\le i<j<k\le n} \mu_{ijk}(M,\sigma).e_i\we e_j\we e_k$.\\ 
The following lemma is a direct consequence of computations contained in the preceding proof (Case 4). 
\begin{lem} \label{lem:miln}
The following diagram commutes 
\begin{diagram} 
\mathcal{A}_1(P_n) & \rOnto^{\varphi_1} & \ov{\mathcal{SL}}^{hb}_1(n) \\  
& \rdOnto & \dTo>{\mu_3} \\ 
& & \Lambda^3H. 
\end{diagram} 
\end{lem}
Lemmas \ref{lem:miln} and \ref{lem:slam} can then be summarized as follows.
\begin{prop} \label{thY2} 
The diagram  
\begin{diagram}  
\mathcal{A}_1(P_n) & \rTo^{\varphi_1} & \ov{\mathcal{SL}}_1^{hb}(n) \\
 & \rdTo<{\rho}  & \dTo>{(\mu_3,\tau)} \\  
 & & \Lambda^{3} H \oplus V
\end{diagram}
commutes, and all of its arrows are \emph{isomorphisms}.
\end{prop} 
More precisely, Lem. \ref{lem:miln} and \ref{lem:slam} imply the commutativity. The fact that $\varphi_1$ (and thus $(\mu_3,\tau)$) 
is an isomorphism follows. 

We are now ready to prove Theorem \ref{corY2}. Assertion (c)$\Longrightarrow$(a) is indeed a direct consequence of Prop. \ref{thY2}. 
As outlined in Rem. \ref{remthm}, assertion (a)$\Longrightarrow$(b) is a general fact, which follows from the definition of a 
finite type invariant. Let us  prove that (b) implies (c) by showing that in fact any homomorphism of Abelian groups 
$\ov{\mathcal{SL}}^{(hb)}_1(n) \rTo^f A$ gives a degree $1$ invariant. 
Let $(M,\sigma)$ be a $n$-string link in a homology ball and let $G_1,G_2$ be some disjoint $Y$-graphs for $(M,\sigma)$.
We aim to show that:
\begin{equation}\label{eq:deg1}
f(M,\sigma)-f\left((M,\sigma)_{G_1}\right)-f\left((M,\sigma)_{G_2}\right)+f\left((M,\sigma)_{G_1\cup G_2}\right)=0.
\end{equation}
Let $G$ be a collection of disjoint $Y$-graphs for $(1_{D^2},1_n)$ such that $(M,\sigma)=\left(1_{D^2},1_n \right)_G$ 
(up to $Y_2$-equivalence). By possibly isotoping $G_1$ and $G_2$ in $M\setminus \sigma$, they are disjoint from $G$.
We then put $(M_i,\sigma_i)=\left( (1_{D^2},1_n) \right)_{G_i}$. Up to $Y_2$-equivalence, $(M,\sigma)_{G_i}=(M,\sigma)\cdot (M_i,\sigma_i)$ 
and $(M,\sigma)_{G_1\cup G_2}=(M,\sigma)\cdot (M_1,\sigma_1) \cdot (M_2,\sigma_2)$. Equation (\ref{eq:deg1}) follows then from the additivity of $f$.
%
%
\section{On the Milnor-Johnson correspondence} \label{slam} %
%
In this section, we study the relation between the  Goussarov-Habiro theory for framed string links in homology balls 
and  this theory for homology cylinders. 
Let us start with a short reminder on the latter. 
\subsection{Homology cylinders} \label{cyl}
Let $\Sigma_{g,1}$ be a compact connected oriented surface of genus $g$ with $1$ boundary component. 

A \emph{homology cylinder} $M$ over $\Sigma_{g,1}$ is a  homology cobordism with an extra homological triviality  
condition \cite{GL,habi,L}. Alternatively, it can be defined as follows: a homology cylinder $M$ over $\Sigma_{g,1}$ 
is a 3-manifold obtained from $\Sigma_{g,1}\times I$ by surgery along some claspers, that is, $M\sim_{Y_1} \Sigma_{g,1}\times I$. \\
The set of homology cylinders up to orientation-preserving diffeomorphisms is denoted by $\mathcal{HC}(\Sigma_{g,1})$. 
It is equipped with a structure of monoid, with product given by the stacking product and  
with $\Sigma_{g,1}\times I$ as unit element.

There is a descending filtration of monoids 
$$\mathcal{HC}(\Sigma_{g,1})=\mathcal{C}_1(\Sigma_{g,1})\supset  \mathcal{C}_2(\Sigma_{g,1})\supset \cdots \supset 
\mathcal{C}_k(\Sigma_{g,1})\supset \cdots$$
where $\mathcal{C}_k(\Sigma_{g,1})$ is the submonoid of all homology cylinders which are $Y_k$-equivalent to $1_{\Sigma_{g,1}}$. 
Moreover, as in the string link case, the quotient monoid  
$\overline{\mathcal{C}}_k(\Sigma_{g,1}):=\mathcal{C}_k(\Sigma_{g,1})/Y_{k+1}$ 
is an Abelian group for every $k\geq 1$. 

As mentioned in \cite{GL,habi}, the Torelli group $\mathcal{T}_{g,1}$ of $\Sigma_{g,1}$ (the isotopy classes of self-diffeomorphisms of 
$\Sigma_{g,1}$ inducing an isomorphism in homology) naturally imbeds in $\mathcal{HC}(\Sigma_{g,1})$ via the mapping cylinder construction, 
and we can extend classical applications on the Torelli group to the realms of homology cylinders. 
In particular, we can extend the \emph{first Johnson homomorphism} $\eta_{1}$ and the \emph{Birman-Craggs homomorphism} $\beta$, originally used 
by D. Johnson in \cite{J1,J2} for the computation of the the Abelianized Torelli group. 
In \cite{MM}, it is shown that these extensions actually are the degree $1$ Goussarov-Habiro finite type invariants for homology 
cylinders. 
%
\begin{thm}[\cite{MM}] \label{hcY2}
Let $M$ and $M'$ be two homology cylinders over $\Sigma_{g,1}$. The following assertions are equivalent:
\begin{enumerate}
\item[(a)] $M$ and $M'$ are  $Y_2$-equivalent;
\item[(b)] $M$ and $M'$ are not distinguished by degree $1$ Goussarov-Habiro finite type invariants;
\item[(c)] $M$ and $M'$ are not distinguished by the first Johnson homomorphism nor the Birman-Craggs homomorphism.
\end{enumerate}
\end{thm} 
This is proved, as in \S \ref{y2slbh}, by computing the abelian group $\ov{\mathcal{C}}_1(\Sigma_{g,1})$ in a graphical way.
More precisely, the authors define (in a strictly similar way) a space of diagrams $\mathcal{A}_1(P_{g,1})$ and a surjective 
surgery map  
$\mathcal{A}_1(P_{g,1}) \rTo^{\psi_1} \ov{\mathcal{C}}_1(\Sigma_{g,1})$, 
which actually is an isomorphism, with inverse given by $\eta_{1}$ and $\beta$. 
\subsection{From homology cylinders to string links}
This result on homology cylinders over $\Sigma_{g,1}$ looks quite similar to Thm. \ref{corY2} on framed $n$-string links 
in homology balls, and suggests a strong analogy between these objects.\\ 
This correspondence homology cylinders/string links has been studied by N. Habegger \cite{habe}: \emph{via} a certain 
geometric construction relating these objects, Johnson homomorphisms coincides with Milnor's numbers.
This result is refered to as the \emph{Milnor-Johnson correspondence}. 
More precisely, Habegger shows that there exists a bijection between the sets 
$\mathcal{HC}(\Sigma_{g,1})$ and $\mathcal{SL}_1^{hb}(2g)$ which produces an isomorphism of Abelian groups 
$$b: \ov{\mathcal{C}}_1(\Sigma_{g,1})\rTo^{\simeq} \ov{\mathcal{SL}}_1^{hb}(2g)$$
\noindent such that the Johnson homomorphism $\eta_1$ corresponds to Milnors invariant $\mu_3$ trough $b$. 
%
%
%
Proposition \ref{thY2} allows us to go a bit further.
%
\begin{thm} \label{thslam} 
The homomorphism $\tau$ of Proposition \ref{thY2}, given by the Milnor, Sato-Levine, Arf and Rochlin invariants, is the analogue of the Birman-Craggs 
homomorphism $\beta$ for the Milnor-Johnson correspondence. \\
In other words, $\beta$ and $\tau$ correspond through the isomorphism $b$. 
\end{thm} 

\noindent The proof is given in the next subsection. Actually, we will also give an alternative proof for 
(part of) Habegger's result, based on the theory of claspers. 
\subsection{Birman-Craggs homomorphism for string links: proof of Theorem \ref{thslam}}
\begin{figure}[!h]
\begin{center}
\includegraphics{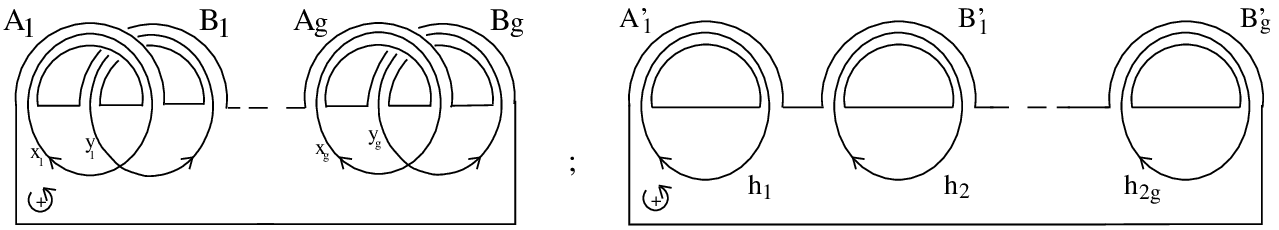}
\caption{} \label{anse}
\end{center}
\end{figure}

Let us recall from \cite{habe} the construction on which the Milnor-Johnson correpondence lies.
Consider the handle decomposition $A_1,B_1,...,A_g,B_g$ of $\Sigma_{g,1}$ as in the left part of Fig. \ref{anse}. 
Likewise, for the $2g$-punctured disk $D^2_{2g}\cong\Sigma_{0,2g+1}$, consider the handle decomposition $\{A'_i,B'_i\}_{i=1}^g$ 
given in the right part of the figure. 

We identify $\Sigma_{g,1}\times I$ with $\Sigma_{0,2g+1}\times I$ using the diffeomorphism $F$ defined by the $g$ isotopies exchanging, in 
$\Sigma_{g,1}\times I$, the second attaching region of the handle $A_i\times I$ and the first attaching region of the handle $B_i\times I$.\\
Now, the product $\Sigma_{0,2g+1}\times I$ can be thought of as (the closure of) the complementary of the $0$-framed trivial $2g$-string link $1_{2g}$ in 
$D^2\times I$.
This defines a bijection between the sets $\mathcal{C}_1(\Sigma_{g,1})$ and $\mathcal{SL}^{hb}_1(2g)$. \\
Indeed, let $G$ be a degree $1$ clasper for $\Sigma_{g,1}\times I$: the pair $(\Sigma_{g,1}\times I;G)$ defines an element of $\mathcal{C}_1(\Sigma_{g,1})$. 
By applying $F$ to this pair, we obtain a clasper $G'$ \emph{of the same degree} for $(\Sigma_{0,2g+1}\times I)\cong 1_{D^2}\setminus 1_{2g}$: the 
triple $\left((1_{D^2},1_n);G' \right)$ defines an element of $\mathcal{SL}^{hb}_1(2g)$. \\
Moreover, though this bijection is not a homomorphism, it produces an isomorphism of Abelian groups
\begin{diagram} 
\ov{\mathcal{C}}_1(\Sigma_{g,1}) & \rTo^{b} & \ov{\mathcal{SL}}^{hb}_1(2g). 
\end{diagram}
\noindent This follows from the following observation. 
Let $M_i$ ($i$=1,2) be an element of $\ov{\mathcal{C}}_1(\Sigma_{g,1})$ obtained from $\Sigma_{g,1}\times I$ by surgery on the degree $1$ clasper $G_i$. 
The product $M_1\cdot M_2$ is mapped by $b$ to an element which is obtained from $(1_{D^2},1_{2g})$ by surgery on the union $G'_1\cup G'_2$, where 
$G'_i$ is the image of $G_i$ under the diffeomorphism $F$ (in particular, $deg(G'_i)=1$). 
Up to $Y_2$-equivalence, we can suppose that these two claspers lie in \emph{disjoint portions of the product $D^2\times I$} ; 
it follows that  
\bc
$(1_{D^2},1_{2g})_{G'_1\cup G'_2}\sim_{Y_2} (1_{D^2},1_{2g})_{G'_1}\cdot (1_{D^2},1_{2g})_{G'_2}=b(M_1)\cdot b(M_2)$.
\ec
Similar arguments show that we actually have an isomorphism of Abelian groups 
$\ov{\mathcal{C}}_k(\Sigma_{g,1})\simeq \ov{\mathcal{SL}}^{hb}_k(2g)$, $\forall$ $k\ge 1$.

At the level of homology, there is an obvious isomorphism between $H_1(\Sigma_{g,1};\mathbf{Z})$ and $H_1(\Sigma_{0,2g+1};\mathbf{Z})$ 
induced by the diffeomorphism $F$. We denote by $H$ these homology groups. 
This isomorphism allows to identify the diagram spaces $\mathcal{A}_1(P_{g,1})$ and $\mathcal{A}_1(P_{2g})$. 
We thus have a commutative diagram
\begin{diagram}
 & & & \mathcal{A}_1(P_{g,1}) & \rTo & \mathcal{A}_1(P_{2g}) & & & \\  
 & & & \dTo<{\psi_1} & & \dTo>{\varphi_1} & & & (D)\\ 
 & & & \ov{\mathcal{C}}_1(\Sigma_{g,1}) & \rTo^{b} & \ov{\mathcal{SL}}^{hb}_1(2g), & & & 
\end{diagram}
whose arrows are isomorphisms. 

Following Notations \ref{notations}, set $H_{(2)}=H\otimes \mathbf{Z}_2$, and $V=\Lambda^2H_{(2)}\oplus H_{(2)}\oplus \mathbf{Z}_2$. 
By considering the inverse maps (in the sense of \cite[Thm. 1.4]{MM} and Prop. \ref{thY2}) of the vertical arrows of $(D)$, we 
easily deduce the following commutative diagram
\begin{diagram} 
\ov{\mathcal{C}}_1(\Sigma_{g,1}) & \rTo^{b} & \ov{\mathcal{SL}}^{hb}_1(2g) \\  
 & \rdTo<{(\eta_1,\beta)}>{\simeq} & \dTo>{(\mu_3,\tau)}<{\simeq} \\ 
 & & \Lambda^3H\oplus V, 
\end{diagram}
\noindent which shows that, \emph{via} the isomorphism $b$, degree $1$ invariants for homology cylinders over $\Sigma_{g,1}$ correspond to those of 
$2g$-string links in homology balls. 
More precisely, we deduce from diagram $(D)$ the following result. 
\begin{lem} \label{lemMJ}
The two following diagrams commute.
\begin{diagram}
\ov{\mathcal{C}}_1(\Sigma_{g,1}) & \rTo^{b} & \ov{\mathcal{SL}}^{hb}_1(2g) & & 
 & \ov{\mathcal{C}}_1(\Sigma_{g,1}) & \rTo^{b} & \ov{\mathcal{SL}}^{hb}_1(2g)\\
 &\rdTo_{\eta_1} & \dTo>{\mu_{3}} & & ;  
 & &\rdTo_{\beta} & \dTo>{\tau} \\
 &  & \Lambda^{3} H & & 
 & & & \Lambda^{3} H_{(2)}\oplus V.\\
\end{diagram} 
\end{lem}
The first diagram recovers Habegger's Milnor-Johnson correspondence (at the lowest level). The second one proves Thm. \ref{thslam}.\\
\noindent \textit{Proof of Lemma \ref{lemMJ}:} \\
Consider in diagram $(D)$ the projections $p : \mathcal{A}_1(P) \rOnto \Lambda^3H$, on the one hand, and 
$T : \mathcal{A}_1(P) \rOnto \Lambda^{3} H_{(2)}\oplus V$ on the other hand, where $\mathcal{A}_1(P)$ denotes either 
$\mathcal{A}_1(P_{g,1})$ or $\mathcal{A}_1(P_{2g})$.  Recall from \cite[Lemma 4.22]{MM} that the diagram 
\begin{diagram} 
\mathcal{A}_1(P_{g,1}) & \rOnto^{\psi_1} & \ov{\mathcal{C}}_1(\Sigma_{g,1}) \\  
& \rdOnto<{p} & \dTo>{\eta_1} \\ 
& & \Lambda^3H
\end{diagram} 
is commutative. This, together with Lemma \ref{lem:miln}, implies the commutativity of the first diagram. 
The second half of the result follows similarly from \cite[Lemma 4.23]{MM} and Lemma \ref{lem:slam}. 
\begin{flushright}
$\square$
\end{flushright}
%
\section{Comparing Goussarov-Habiro and Vassiliev theories} \label{ghv}  %
\begin{figure}[!h]
\begin{center}
\includegraphics{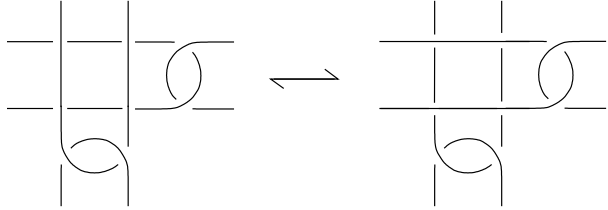}
\caption{A clasp-pass move.} \label{cpm}
\end{center}
\end{figure}

%
For several reasons, it is tempting to compare the results of \S \ref{GHSL} with Vassiliev theory. 
First, as seen in \S \ref{pre}, both Goussarov-Habiro and Vassiliev theories can be defined using claspers 
(with some slight differences). Second, some results in the literature on Vassiliev invariants have 
strong similarities with Theorem \ref{corY2}, namely K. Taniyama and A. Yasuhara's characterization of 
clasp-pass equivalence for algebraically split links in the 3-sphere \cite{TY}, and its analogue for 
string links \cite{moi}. \\
Recall that the clasp-pass equivalence is the equivalence relation on links generated by isotopies and \emph{clasp-pass moves}, 
which are local moves as illustrated in Fig. \ref{cpm}.
As outlined in Rem. \ref{Ckmoves}, the clasp-pass equivalence is actually the same as $C_3$-equivalence, 
which implies $Y_2$-equivalence. 
\subsection{Goussarov-Habiro and Vassiliev invariants of string links}
Let us first consider the string link case. 
Recall that the Casson knot invariant $\varphi(K)$ of a knot $K$ is defined as the $z^2$ coefficient of the Alexander-Conway polynomial 
of $K$, and that its reduction modulo 2 coincides with the Arf invariant $\alpha$ studied in \S \ref{ARF}. \\
Recall also from \cite{moi} the definition of the 2-string link invariant $V_2$. 
Let $\sigma=\sigma_1\cup \sigma_2$ be a 2-string link. Then 
$$
V_2(\sigma) : = \varphi\left(p(\sigma)\right) - \varphi(\sigma_1) - \varphi(\sigma_2), 
$$
where $p(\sigma)$ denotes the plat-closure of $\sigma$: it is the knot obtained by identifying the upper (resp. lower) endpoints of 
$\sigma_1$ and $\sigma_2$. Clearly, $V_2$ is a $\mathbf{Z}$-valued Vassiliev invariant of degree two.\\
We want to relate Thm. \ref{corY2} to the following:
\begin{thm}[\cite{moi}] \label{C3sl}
Let $\sigma$ and $\sigma'$ be two $n$-component algebraically split string links in $D^2\times I$ 
(that is, with all linking numbers zero). Then, the following assertions are equivalent:
\begin{enumerate}
\item[(a)] $\sigma$ and $\sigma'$ are clasp-pass equivalent;
\item[(b)] $\sigma$ and $\sigma'$ are not distinguished by degree $2$ Vassiliev invariants;
\item[(c)] $\sigma$ and $\sigma'$ are not distinguished by Milnor's triple linking numbers, nor the invariant $V_2$ and the Casson knot invariant.
\end{enumerate}
\end{thm} 

We denote by $SL(n)$ the monoid of $n$-string links in $D^2\times I$ up to isotopy (with fixed endpoints), and 
by $SL^{as}(n)$ the submonoid of algebraically split $n$-string links.
When considered up to $C_3$-equivalence, the elements of $SL^{as}(n)$ form an Abelian group, denoted by $\ov{SL}^{as}(n)$. 
\begin{thm} \label{bonus}
The Abelian group $\ov{SL}^{as}(n)$ is surjectively mapped onto the subgroup 
$\ov{\mathcal{SL}}^{(0)}_1(n)\subset \ov{\mathcal{SL}}_1(n)$ of string links in homology balls having vanishing Rochlin's $\mu$-invariant. 
\end{thm}

\emph{Proof: }
Recall from \cite{moi} the isomorphism 
$$ (\mu_3,V_2,\varphi): \ov{SL}^{as}(n)\rTo^{\simeq} \Lambda^3H\oplus S^2H $$
given by the formula
$$ \sum_{1\le i<j<k\le n} \mu_{ijk}(\sigma).e_i\we e_j\we e_k 
+ \sum_{1\le i<j\le n} V_2(\sigma_i\cup \sigma_j).e_i\otimes e_j + \sum_{1\le i\le n} \varphi(\sigma_i).e_i.$$
Here, $S^2H$ is the degree two part of the symmetric algebra of $H$ (we still make use of Notations \ref{notations}).

On the other hand, we saw in \S \ref{y2slbh} the isomorphism of Abelian groups 
$$ (\mu_3,\tau): \ov{\mathcal{SL}}^{hb}_1(n)\rTo^{\simeq} \Lambda^{3} H \oplus \Lambda^{2} H_{(2)}\oplus H_{(2)}\oplus \mathbf{Z}_{2},$$ 
where the $\mathbf{Z}_2$ part is detected by Rochlin's $\mu$-invariant. We thus have the decomposition 
$\ov{\mathcal{SL}}^{hb}_1(n) = \ov{\mathcal{SL}}^{(0)}_1(n)\cup \ov{\mathcal{SL}}^{(1)}_1(n)$,
where $\ov{\mathcal{SL}}^{(\epsilon)}_1(n)$ ($\epsilon=0,1$) is the subset of $\ov{\mathcal{SL}}^{hb}_1(n)$ consisting of 
elements $(M,\sigma)$ such that $R(M)=\epsilon$. \\
In particular, $\ov{\mathcal{SL}}^{(0)}_1(n)$ is an Abelian subgroup and we clearly have an isomorphism
$$(\mu_3,\beta^{(2)},\alpha): \ov{\mathcal{SL}}^{(0)}_1(n)\rTo^{\simeq}  \Lambda^{3} H \oplus \Lambda^{2} H_{(2)}\oplus H_{(2)},$$
given by the formula
$$ \sum_{1\le i<j<k\le n} \mu_{ijk}(M,\sigma).e_i\we e_j\we e_k  + \sum_{1\le i<j\le n} \beta^{(2)}_{ij}(M,\sigma).\ov{e_i}\we \ov{e_j}+ 
\sum_{1\le i\le n} a_{i}(M,\sigma).\ov{e_i}.$$

Now, recall that the $C_2$-equivalence is the same as the $\Delta$-equivalence: as in the link case \cite{MN}, a $n$-string link $\sigma$ is 
$C_2$-equivalent to $1_n$ if and only if it has vanishing linking numbers. So $\ov{SL}^{as}(n)$ is just the set of $C_3$-equivalence 
classes of $n$-string links which are $C_2$-equivalent to $1_n$: given a generator $\sigma$ of $\ov{SL}^{as}(n)$, there is a connected $C$-degree 2 clasper $G_{\sigma}$ 
for $1_n\in 1_{D^2}$ such that $\sigma=(1_n)_{G_{\sigma}}$.
We define a map
\begin{diagram}
\ov{SL}^{as}(n) & \rTo^{\eta} & \ov{\mathcal{SL}}^{(0)}_1(n) 
\end{diagram} 
which consists in \emph{puncturing} 
each disk-leaf of $G_{\sigma}$, that is removing a small disk $d$ such that $1_n$ intersects the disk-leaf at the interior of $d$ ; 
further, equip $1_n$ with $0$-framing.
\begin{figure}[!h]
\begin{center}
\includegraphics{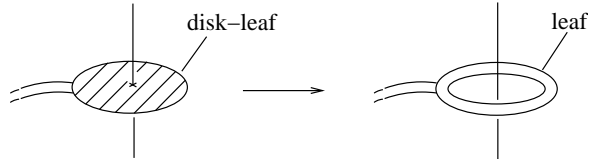}
\caption{The $\eta$ map.} \label{trou}
\end{center}
\end{figure}
As Fig. \ref{trou} shows, puncturing a disk-leaf of $G_{\sigma}$ produces a leaf. $G_{\sigma}$ becomes a Y-graph $\tilde{G}_{\sigma}$, 
and
\bc
$T(\sigma) := (1_{D^2},1_n)_{\tilde{G}_{\sigma}}.\quad $
\ec
\noindent Note that $\eta$ has a non-trivial kernel ; an example is given in Fig. \ref{mec}.
\begin{figure}
\bc
\includegraphics{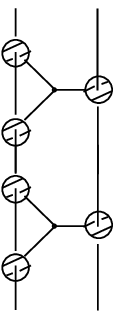}
\caption{An element of $Ker(\eta)$ for $n=2$.} \label{mec}
\ec
\end{figure}
\noindent It follows from the proofs of Thm. \ref{corY2} and \ref{C3sl} that we have a commutative diagram
\begin{diagram}
\ov{SL}^{as}(n) & \rTo^{\eta} & \overline{\mathcal{SL}}^{(0)}_1(n)  \\
\dTo<{(\mu_3,V_2,\varphi)}>{\simeq} &  & \dTo<{\simeq}>{(\mu_3,\beta^{(2)},\alpha)} \\
\Lambda^{3}H\oplus S^{2}H & \rOnto^{t} & \Lambda^{3} H \oplus\Lambda^{2} H_{(2)}\oplus H_{(2)}, 
\end{diagram}
where $f$ is the surjective map given by the identity on $\Lambda^{3}H$, and by 
$$ f(e_i\otimes e_j)=\ov{e_i}\we \ov{e_j} \textrm{ if $i\ne j$, and } \\ 
f(e_i\otimes e_i)=\ov{e_i} \textrm{ otherwise } $$
on $S^{2}H$. It follows that $\eta$ is also surjective. $\square$

Moreover, the maps $(\mu_{3},V_2,c_2)$ and $(\mu_{3},\beta^{(2)},a)$ 
coincide \emph{via} the surjective map $\eta$ (and $t$). In particular, it follows that 
$$V_2 \equiv \beta \textrm{ (mod $2$).}$$
\noindent However, these invariants are distinct over $\mathbf{Z}$, as mentionned in \cite[Rem 2.7]{moi}.
%
%
%
%
\subsection{The case of links} \label{link}	
In the case of links, we know the following on clasp-pass equivalence.
\begin{thm}[\cite{TY}, Thm. 1.4] \label{C3l}
Let $L$ and $L'$ be two $n$-component algebraically split links in $S^3$. 
The following assertions are equivalent:
\begin{enumerate}
\item[(a)] $L$ and $L'$ are clasp-pass equivalent;
\item[(b)] $L$ and $L'$ are not distinguished by Milnor's triple linking numbers, nor the mod 2 reduction of 
the Sato-Levine invariant and the Casson knot invariant.
\end{enumerate}
\end{thm} 
\noindent As for $Y_2$-equivalence, one can check (using Thm \ref{corY2} and its proof) the following corollary, 
characterizing $Y_2$-equivalence for algebraically split links in homology spheres.
\begin{cor}
Let $(M,L)$ and $(M',L')$ be two $n$-component algebraically split links in homology spheres.  Then, the following assertions are equivalent:
\begin{enumerate}
\item[(a)] $(M,L)$ and $(M',L')$ are $Y_2$-equivalent;
\item[(b)] $(M,L)$ and $(M',L')$ are not distinguished by Milnor's triple linking numbers, nor the mod 2 reduction of 
the Sato-Levine invariant, the Arf invariant and Rochlin's $\mu$-invariant.
\end{enumerate}
\end{cor}
This result is related to Thm. \ref{C3l} in a similar way as Thm. \ref{corY2} is related to Thm. \ref{C3sl}.
However, unlike in the string link case, there is no natural group or monoid structure on the sets of $C_k$ or $Y_k$-equivalence classes of links. 
%
%
%
%
%
%
\appendix
\section{Tubing Seifert surfaces.} \label{appendix}
Let us consider the 2-component link $L=L_1\cup L_2$ in a genus $4$ handlebody $N$ depicted in Fig. \ref{fig:link}. 
We fix an orientation on $N$ and embed it in $S^3$. 
Let $K$ be an oriented knot in $S^3$ disjoint from $N$, and let $S$ be a Seifert surface for $K$: in general, 
$S$ may intersect $K$. In this appendix we explain the general procedure to construct, starting from $S$, a new Seifert surface 
for $K$ which is disjoint from $L$. 

First, we fix some more notations. The handlebody $N$ can be regarded as a ball $B$ with $4$ handles $D^2\times I$ attached. 
The two handles intersecting $L_1$ are denoted by $H_1$ and $H_2$, and we denote by $H_3$ and $H_4$ the other two ; the handles 
are numbered clockwise in Fig. \ref{fig:link}, so that $H_1$ is in the lower left corner of the figure. 
Up to isotopy, we can suppose that $S$ is disjoint from $B$, that is $S$ only intersects $N$ at its handles, along copies of 
$D^2\times\{t\}$ ; $t\in I$. When the orientation of $S$ is compatible with the orientation of $N$ along the intersection disk, 
we call it a \emph{positive intersection}. Otherwise, we call it a \emph{negative intersection}. For $1\le i\le 4$, we denote 
respectively by $p_i$ and $n_i$ the number of positive and negative intersections between the surface $S$ and the handle $H_i$.\\ 
In view of the symmetry of the link $L$, we only have to deal with (say) the handles $H_1$ and $H_2$ (the handles $H_3$ and $H_4$ 
can be treated independently, in a similar way). 

\begin{figure}[!h]
\begin{center}
\includegraphics{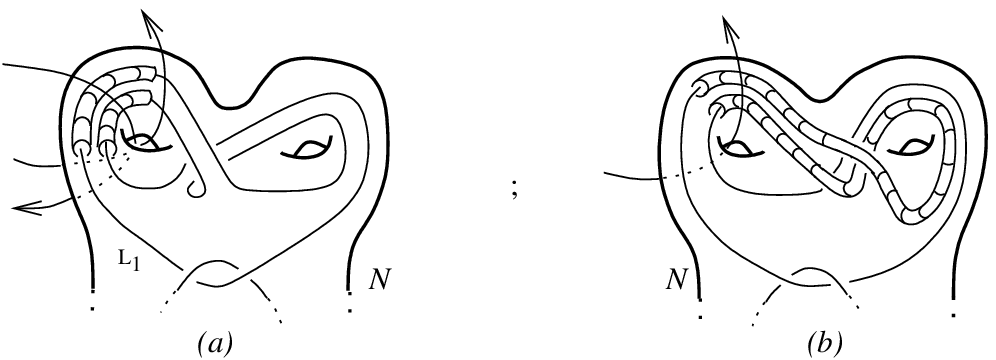}
\caption{} \label{tube}
\end{center}
\end{figure}

First, observe that if $S$ intersects $H_1$ twice, with the opposite orientation, we can add two tubes to $S$ as shown in Fig. 
\ref{tube} (a), so that the new Seifert surface $\tilde{S}$ satisfies $|L\cap S|=|L\cap \tilde{S}|+4$.

\noindent Likewise, we can always add $|p_i - n_i|$ such pairs of tubes to $S$ in $H_i$ ($i=1,2$), by eventualy nesting them, 
so that in each handle the remaining intersections all have the same sign. 
So we can suppose that $p_1.n_1=p_2.n_2=0$. 
Suppose further that $n_1=0$ (the case $p_1=0$ is equivalent, due to the symmetry of $L_1$). \\
If $p_1=p_2=n_2=0$, we are done. Otherwise, there are essentially $4$ different cases to study. 

1. Suppose that $p_2=n_2=0$. 
In this case $S$ is disjoint from the handle $H_1$. We can thus remove all the elements of $S\cap L_1$ by successively 
attaching and nesting $p_1$ tubes as depicted in Fig. \ref{tube}(b). These tubes will be called tubes \emph{of type 1}. 
When the two attaching circles of a tube $t$ of type 1 lie in a disk $D^2\times\{t\}$ of the handle $H_i$, we simply say that $t$ 
is \emph{attached in} $H_i$ ($i=1,2$).

2. Suppose that $p_1=0$. 
This case is equivalent to the first one: $S$ is disjoint from the handle $H_2$, so we can freely attach $p_2+n_2$ tubes 
of type 1 in $H_1$.

\begin{figure}[!h]
\begin{center}
\includegraphics{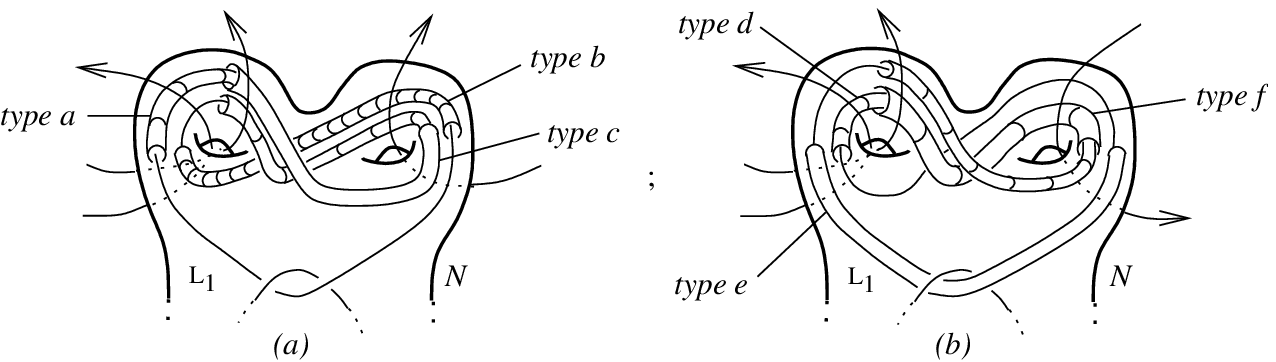}
\caption{} \label{tubes}
\end{center}
\end{figure}

3. Suppose that $p_1$ and $p_2$ are non-zero.
In this case, $S$ always intersects $N$ with the same sign. Fig. \ref{tubes} $(a)$ illustrates the case $p_1=2$ and $p_2=1$.\\
In general, we attach in a similar way $p$ tubes of type $a$, $m$ tubes of type $b$ and $m$ tubes of type $c$ (following the notations 
of the figure), where $p:=|p_1-p_2|$ and $m:=max(p_1,p_2)-p$.

4. Suppose that $p_1$ and $n_2$ are non-zero. 
Fig. \ref{tubes} $(b)$ illustrates the case $p_1=2$ and $n_2=1$. As for the previous case, we deal with the general 
situation by attaching and nesting the same three types of tubes. Namely, we attach $|p_1-n_2|$ tubes of type $d$ and 
$\left(max(p_1,n_2)-|p_1-n_2| \right)$ tubes of type $e$ and $f$.

The obtained surface is the required new Seifert surface for $K$.
%
%
%
%
%
%

\vspace{0.2cm}
\noindent
\footnotesize{Commutative diagrams were drawn with Paul Taylor's package.}
\end{document}